\documentclass[11pt]{article}
\usepackage{t1enc}
\usepackage{lmodern}
\usepackage[T2A,T1]{fontenc}
\usepackage[utf8]{inputenc}
\usepackage[russian,english]{babel}
\usepackage{amsmath,amssymb,amsfonts,amsthm,mathrsfs,textcomp,url,bbm,cancel,comment,enumerate}
\usepackage{graphicx}
\usepackage{float}
\usepackage[all,ps]{xy}
\usepackage[colorlinks,urlcolor=cyan,citecolor=blue,linkcolor=blue]{hyperref}
\usepackage[dvipsnames]{xcolor}
\usepackage[margin=3cm]{geometry}
\usepackage{longtable}
\usepackage{tikz}
\makeatletter
\newcommand{\affil}[1]{\gdef\@affil{\textsc{#1}}}
\newcommand{\address}[1]{\gdef\@address{#1}}
\newcommand{\email}[1]{\gdef\@email{\url{#1}}}
\newcommand{\@endstuff}{\par\vspace{\baselineskip}\noindent\footnotesize
  \begin{tabular}{@{}l}
    \@affil\\
    \@address\\
    \textit{E-mail address:} \@email
  \end{tabular}}
\AtEndDocument{\@endstuff}
\makeatother

\usetikzlibrary{matrix,calc,decorations,positioning}
\pgfkeys{/tikz/.cd,
    alt double distance/.initial=5pt,
    alt double step/.initial=1pt,}
\pgfdeclaredecoration{double deco}{initial}
{
\state{initial}[width=\pgfkeysvalueof{/tikz/alt double step},next state=cont] {
    \pgfmoveto{\pgfpoint{\pgfkeysvalueof{/tikz/alt double step}}{\pgfkeysvalueof{/tikz/alt double distance}/2}}
    \pgfpathlineto{\pgfpoint{0.3\pgflinewidth}{\pgfkeysvalueof{/tikz/alt double distance}/2}}
    \pgfpathmoveto{\pgfpoint{0.3\pgflinewidth}{-\pgfkeysvalueof{/tikz/alt double distance}/2}}
    \pgfpathlineto{\pgfpoint{1pt}{-\pgfkeysvalueof{/tikz/alt double distance}/2}}
    \pgfcoordinate{lastup}{\pgfpoint{1pt}{\pgfkeysvalueof{/tikz/alt double distance}/2}}
    \pgfcoordinate{lastdown}{\pgfpoint{1pt}{-\pgfkeysvalueof{/tikz/alt double distance}/2}}}
  \state{cont}[width=\pgfkeysvalueof{/tikz/alt double step}]{
     \pgfmoveto{\pgfpointanchor{lastup}{center}}
     \pgfpathlineto{\pgfpoint{\pgfkeysvalueof{/tikz/alt double step}}{\pgfkeysvalueof{/tikz/alt double distance}/2}}
     \pgfcoordinate{lastup}{\pgfpoint{\pgfkeysvalueof{/tikz/alt double step}}{\pgfkeysvalueof{/tikz/alt double distance}/2}}
     \pgfmoveto{\pgfpointanchor{lastdown}{center}}
     \pgfpathlineto{\pgfpoint{\pgfkeysvalueof{/tikz/alt double step}}{-\pgfkeysvalueof{/tikz/alt double distance}/2}}
     \pgfcoordinate{lastdown}{\pgfpoint{\pgfkeysvalueof{/tikz/alt double step}}{-\pgfkeysvalueof{/tikz/alt double distance}/2}}}
  \state{final}[width=0pt]
  { 
    \pgfmoveto{\pgfpointdecoratedpathlast}}}
\tikzset{alt double/.style={decorate,decoration=double deco}}
\allowdisplaybreaks

\newcommand{\toverset}[2]{%
\mathop{#2}\limits^{\vbox to -.1ex{\kern-0.4ex\hbox{$\scriptstyle #1$}\vss}}}
\newcommand{\tightoverset}[2]{%
\mathop{#2}\limits^{\vbox to -.5ex{\kern-0.4ex\hbox{$\scriptstyle #1$}\vss}}}
\renewcommand{\underset}[2]{%
\mathop{#2}\limits_{\vbox to -.5ex{\kern-1.6ex\hbox{$\scriptstyle #1$}\vss}}}
\newcommand{\tightunderset}[2]{%
\mathop{#2}\limits_{\vbox to -.5ex{\kern-1.8ex\hbox{$\scriptstyle #1$}\vss}}}
\newcommand{\xra}[1]{%
\mathop{\xrightarrow{~#1~}}}

\newcommand{\congc}[1]{%
\mathop{\tightunderset{\mathscr{C}_{#1}}{\cong}}}

\def\Z{\mathbb{Z}}

\def\R{\mathbb{R}}

\def\RP{\mathbb{R}P}
\def\CP{\mathbb{C}P}

\def\CC{\mathscr{C}}

\def\NN{\mathfrak{N}}

\def\into{\hookrightarrow}
\def\imto{\looparrowright}

\def\ol{\overline}

\def\la{\langle}
\def\ra{\rangle}

\def\pr{\mathop{\rm pr}}

\def\O{\mathrm{O}}

\def\SO{\mathrm{SO}}

\def\Imm{\textstyle{\mathop{\rm Imm}}}

\def\Mor{\textstyle{\mathop{\rm Mor}}}

\def\Prim{\textstyle{\mathop{\rm Prim}}}

\renewcommand{\Lambda}{\det}
\newenvironment{prf}%
{\par\noindent\textbf{Proof.\enspace\ignorespaces}}%
{~$\square$\par}%
{\par\noindent\textit{Claim.\enspace\ignorespaces}\begin{em}}%
{\end{em}\par}%
{\medskip\par\noindent\textit{Proof.\enspace\ignorespaces}}%
{~$\diamond$\par\medskip}%
x{\medskip\par\noindent\textit{Remark.\enspace\ignorespaces}}%
{\par\medskip}%
{\par\noindent\begin{small}\textbf{Keywords.\enspace\ignorespaces}}%
{\end{small}\par}%
{\par\noindent\begin{small}\textbf{Acknowledgement.\enspace\ignorespaces}}%
{\end{small}\par}%
\newtheorem{thm}{Theorem}[section]%
\newtheorem{lemma}[thm]{Lemma}%
\newtheorem{prop}[thm]{Proposition}%
\theoremstyle{definition}
\newtheorem{defi}[thm]{Definition}%
\newtheorem{rmk}[thm]{Remark}%
\newcommand{\ZeroRoman}[1]{%
\ifcase\value{#1}\relax 0\else\Roman{#1}\fi}

\newcounter{t}

\numberwithin{equation}{section}
\title{Computations on cobordism groups of projected immersions}
\author{András Csépai}
\affil{Eötvös Loránd University, Institute of Mathematics}
\address{Pázmány Péter sétány 1/c, Budapest, H-1117 Hungary}
\email{csepai.andras112358@gmail.com}
\date{}
\begin{document}

\maketitle

\begin{abstract}
  We use a classifying space construction by Szűcs and Terpai \cite{ctrl} to extend Salomonsen's exact sequence from cobordism groups of immersions to cobordism groups of hyperplane projections of immersions with prescribed singularities. We apply this to obtain results on these cobordism groups in small and large codimensional cases.
\end{abstract}

\section{Introduction}
\label{sec:intro}

A \textit{projected immersion}, or shortly \textit{prim map}, is defined as a smooth map $f\colon M\to P$ between manifolds such that there is an immersion $g\colon M\imto P\times\R$ for which $f$ is the composition $\pr_P\circ g$ where $\pr_P$ denotes projection to $P$. Prim maps are maps with singularities and their cobordism groups are natural generalisations of cobordism groups of immersions and fit into the larger theory of cobordisms of singular maps; see e.g. \cite{hosszu}, \cite{ando}, \cite{kal}, \cite{sad}, \cite{ctrl} (and the references therein) for different approaches. 

Throughout this paper $G$ will denote either of the groups $\O$ and $\SO$ (although everything goes through without change if $G$ is any stable group in the sense of \cite[section 8.2]{stabgr}, i.e. a direct limit of linear groups with certain properties which imply that it makes sense to talk about cobordism groups with $G$-structure). We shall denote by $\gamma_{k+1}^G$ the tautological vector bundle over $BG(k+1)$ and by $\zeta^G_r$ the pullback of $\gamma_{k+1}^G$ by the projection of $S(r+1)\gamma_{k+1}^G$, that is, the sphere bundle of $\gamma_{k+1}^G\oplus\ldots\oplus\gamma_{k+1}^G$ with $r+1$ summands. Note that here $r=\infty$ also makes sense: $\zeta^G_\infty$ is simply the bundle $\gamma_{k+1}^G$. The cobordism groups of immersions and prim maps $M^n\to\R^{n+k}$ with normal bundle induced from a vector bundle $\xi$ in the case of immersions, and with normal $G$-structure and at most $\Sigma^{1_r}$-type singularities in the case of prim maps, will be denoted by $\Imm^\xi(n,k)$ and $\Prim^G_r(n,k)$ respectively. Note that if $s<r$, there is a natural forgetful homomorphism $\varphi_s^r\colon\Prim^G_s(n,k)\to\Prim^G_r(n,k)$.

Szűcs and Terpai \cite{ctrl} proved the isomorphism $\Prim^G_r(n,k)\cong\Imm^{\zeta^G_r}(n,k+1)$ for all $r$. Using this, we will prove:

\begin{thm}
  \label{thm:main}
  For all $0\le s<r\le\infty$ the forgetful homomorphism $\varphi_s^r$ fits into a long exact sequence
  \begin{alignat*}2
    \ldots&\to\Prim^G_s(n,k)\xra{\varphi_s^r}\Prim^G_r(n,k)\to\Imm^{(s+2)\zeta^G_{r-s-1}}(n-(s+1)(k+1),(s+2)(k+1))\to\\
          &\to\Prim^G_s(n-1,k)\to\ldots
  \end{alignat*}
\end{thm}

The paper is organised as follows: in section \ref{sec:cob} we recall the definitions of singularities and cobordism groups of immersions and prim maps with prescribed singularities; in section \ref{sec:prf} we prove theorem \ref{thm:main} and make a few remarks on it; then in section \ref{sec:comp} we apply this theorem to complete the computation of the cobordism groups $\Prim^G_r(n,k)$ modulo small torsion groups where $k$ is either at most $1$ or at least $\frac{n-2}2$.

\section{Singularities and cobordism of prim maps}
\label{sec:cob}


\begin{defi}
  For a smooth map $f\colon M\to P$ and a number $r$ we define $\Sigma^r(f)\subset M$ to be the set of points $p\in M$ where the corank of the differential of $f$ (i.e. $\dim\ker df_p$) is $r$. For generic maps $f$ the subset $\Sigma^r(f)\subset M$ is a manifold, hence we can consider for any $s\le r$ the subset $\Sigma^{r,s}(f):=\Sigma^s(f|_{\Sigma^r(f)})$ and iterating this process yields the definition of the subset $\Sigma^S(f)$ for any weakly decreasing sequence $S$ of non-negative integers. This classification of singular points defines the \textit{Thom--Boardman singularity type} $\Sigma^S$.
\end{defi}

\begin{defi}
  A map $f\colon M\to P$ is said to be a \textit{Morin map} if $\dim\ker df_p\le 1$ for any point $p\in M$. In this case, each point $p\in M$ belongs to $\Sigma^{1_r}(f)$ for some $r\ge0$ where $\Sigma^{1_r}$ denotes $\Sigma^{1,\ldots,1,0}$ with $r$ digits $1$. A Morin map $f$ is called a \textit{$\Sigma^{1_r}$-map} if $\Sigma^{1_{r+1}}(f)$ is empty.
\end{defi}

\begin{rmk}
  Morin \cite{mor} proved that for a Morin map $f\colon M\to P$, the germs of $f$ at two points $p,q\in M$ coincide up to local coordinate change precisely if they both belong to $\Sigma^{1_r}(f)$ for some $r$. This local singularity defined by $\Sigma^{1_r}$ is also commonly denoted by $A_r$.
\end{rmk}

Now if $f\colon M\to P$ is a prim map, then it is Morin, moreover, the kernel line bundle $\ker df$ over $\Sigma^1(f)$ is trivial. In this paper we assume prim maps $f$ to be equipped with a fixed trivialisation of $\ker df$ which corresponds to fixing the regular homotopy class of their immersion lifts $g\colon M\imto P\times\R$ (i.e. the maps such that $f=\pr_P\circ g$); we will consider projections of immersions with different trivialisations of the kernel bundle to be different prim maps.

\begin{rmk}
  \label{rmk:prim}
  It is not hard to see that the opposite of the above observation also holds, namely a Morin map with trivialised kernel line bundle is prim (see e.g. \cite{immo}).
\end{rmk}

\begin{defi}
  \label{defi:primcob}
  Fix a number $0\le r\le\infty$. Two prim maps $f_0\colon M_0^n\to P^{n+k}$ and $f_1\colon M_1^n\to P^{n+k}$ with closed source manifolds and at most $\Sigma^{1_r}$-type singularities are said to be \textit{prim $\Sigma^{1_r}$-cobordant} if
  \begin{enumerate}[(i)]
  \item there is a manifold $W^{n+1}$ with boundary $\partial W=M_0\sqcup M_1$,
  \item there is a prim $\Sigma^{1_r}$-map $F\colon W\to P\times[0,1]$ such that for $i=0,1$ the restriction $F|_{M_i}$ is $f_i\colon M_i\to P\times\{i\}$.
  \end{enumerate}
  If the stable normal bundles of $f_0$, $f_1$ and $F$ above are equipped with $G$-structures, then we call them \textit{prim $G$-$\Sigma^{1_r}$-cobordant}. This cobordism is an equivalence relation and we denote the cobordism class represented by a map $f$ by $[f]$ (suppressing $G$ and $r$ in the notation).
\end{defi}

\begin{defi}
  \label{defi:cobgroup}
  For two prim $\Sigma^{1_r}$-maps $f\colon M^n\to\R^{n+k}$ and $g\colon N^n\to\R^{n+k}$ with normal $G$-structures, forming their disjoint union defines an Abelian group operation $[f]+[g]:=[f\sqcup g]$ on the set of prim $G$-$\Sigma^{1_r}$-cobordism classes of maps of $n$-manifolds to $\R^{n+k}$; this group will be denoted by $\Prim_r^G(n,k)$.
\end{defi}

\begin{defi}
  For $0\le s<r\le\infty$, assigning to a prim $G$-$\Sigma^{1_s}$-cobordism class $[f]$ the prim $G$-$\Sigma^{1_r}$-cobordism class of $f$ defines the \textit{forgetful homomorphism} $\varphi_s^r\colon\Prim^G_s(n,k)\to\Prim^G_r(n,k)$.
\end{defi}

The analogue of the cobordism group $\Prim^G_r(n,k)$ for general (non-prim) Morin maps can be defined analogously to definitions \ref{defi:primcob} and \ref{defi:cobgroup}; we will denote this group by $\Mor_r^G(n,k)$.

The special case $r=0$ of definition \ref{defi:primcob} yields the cobordism of immersions with stable normal $G$-structures whose cobordism group is denoted by $\Imm^G(n,k)$. In this case not only the stable, but also the actual normal bundle can be specified; if $r=0$ and the normal bundles of $f_0$, $f_1$ and $F$ in definition \ref{defi:primcob} are induced from a vector bundle $\xi$, we get the definition of cobordism of immersions with normal $\xi$-structures, i.e. $\xi$-immersions. This cobordism group is denoted by $\Imm^\xi(n,k)$. Note that we have $\Imm^G(n,k)=\Imm^{\gamma^G_k}(n,k)$ where $\gamma^G_k$ is the tautological vector bundle over $BG(k)$.

\section{Proof of theorem \ref{thm:main}}
\label{sec:prf}

Recall that for any vector bundle $\xi$, there is a Pontryagin--Thom type isomorphism (due to Wells \cite{cobimm}) $\Imm^\xi(n,k)\cong\pi_{n+k}(\Gamma T\xi)\cong\pi^s_{n+k}(T\xi)$, where $T\xi$ denotes the Thom space of $\xi$ and $\Gamma$ is the functor $\Omega^\infty S^\infty$. In other words, the classifying space of cobordisms of $\xi$-immersions is $\Gamma T\xi$. Similarly, cobordisms of singular maps also have classifying spaces, i.e. spaces $X$ such that the appropriate cobordism group of maps $M^n\to\R^{n+k}$ is isomorphic to $\pi_{n+k}(X)$ (see e.g. \cite{hosszu} or \cite{sad} for general constructions). In \cite{ctrl} Szűcs and Terpai identified the classifying space of the prim cobordism groups $\Prim^G_r(n,k)$ for $0\le r\le\infty$ with $\Omega\Gamma T\zeta^G_r$ yielding
$$\Prim^G_r(n,k)\cong\pi^s_{n+k+1}(T\zeta^G_r)\cong\Imm^{\zeta^G_r}(n,k+1).$$

Recall that $\zeta^G_r$ is defined as the pullback $p^*_{r+1}\gamma^G_{k+1}$ where $p_{r+1}\colon S(r+1)\gamma^G_{k+1}\to BG(k+1)$ is the projection of the sphere bundle; for $r=\infty$ this bundle has contractible fibre, hence then $p_{r+1}$ is a homotopy equivalence. The geometric interpretation of $\zeta^G_r$-immersions $g\colon M^n\imto\R^{n+k+1}$ is that they are immersions with normal $G$-structure, endowed with $r+1$ sections $\sigma_1,\ldots,\sigma_{r+1}$ of $\nu_g$ which have no common zeros. The identification of the cobordism group $\Imm^{\zeta^G_r}(n,k+1)$ with $\Prim^G_r(n,k)$ on the level of representatives of cobordism classes assigns to the tuple $(g,\sigma_1,\ldots,\sigma_{r+1})$ a prim map $f\colon M\to\R^{n+k}$ such that for each $i$ the set $\Sigma^{1_i}(f)$ is $\sigma_1^{-1}(0)\cap\ldots\cap\sigma_i^{-1}(0)$ (hence $\Sigma^{1_{r+1}}(f)$ is empty).

\begin{lemma}
  \label{lemma:main}
  For any $0\le s<r\le\infty$ 
  there is a cofibration
  $$T\zeta^G_s\into T\zeta^G_r\to T(s+2)\zeta^G_{r-s-1}.$$
\end{lemma}

\begin{prf}
  Consider the disk bundle $D\zeta^G_r$ of $\zeta^G_r$ which is the pullback $p_{r+1}^*D\gamma^G_{k+1}$. Now $S(r+1)\gamma^G_{k+1}$ is the union of the disk bundles of the pullbacks of $(s+1)\gamma^G_{k+1}$ and $(r-s)\gamma^G_{k+1}$ over $S(r-s)\gamma^G_{k+1}$ and $S(s+1)\gamma^G_{k+1}$ respectively, that is,
  $$S(r+1)\gamma^G_{k+1}=D(s+1)\zeta^G_{r-s-1}\cup D(r-s)\zeta^G_s$$
  glued along their common sphere bundles. Thus $D\zeta^G_r$ is the union of $D(s+2)\zeta^G_{r-s-1}$ and $D(r-s+1)\zeta^G_s$. We obtain the Thom space $T\zeta^G_r$ by factoring $D\zeta^G_r$ by $S\zeta^G_r$ and if we also factor this Thom space by the subspace $D(r-s+1)\zeta^G_s/S\zeta^G_r|_{D(r-s)\zeta^G_s}$, we get $T(s+2)\zeta^G_{r-s-1}$. But $D(r-s+1)\zeta^G_s$ deformation retracts to $D\zeta^G_s$ and this retraction takes the subspace $S\zeta^G_r|_{D(r-s)\zeta^G_s}$ to $S\zeta^G_s$, hence $T(s+2)\zeta^G_{r-s-1}$ is the factor of $T\zeta^G_r$ by $T\zeta^G_s$. This is what we wanted.
\end{prf}\medskip

The long exact sequence of the cofibration in lemma \ref{lemma:main} in stable homotopy groups is
\begin{alignat*}2
  \ldots&\to\Imm^{\zeta^G_s}(n,k+1)\to\Imm^{\zeta^G_r}(n,k+1)\to\\
        &\to\Imm^{(s+2)\zeta^G_{r-s-1}}(n-(s+1)(k+1),(s+2)(k+1))\to\Imm^{\zeta^G_s}(n-1,k+1)\to\ldots
\end{alignat*}
Here the homomorphism $\Imm^{\zeta^G_s}(n,k+1)\to\Imm^{\zeta^G_r}(n,k+1)$ is induced by the natural inclusion $T\zeta^G_s\subset T\zeta^G_r$, hence it assigns to the cobordism class of an immersion equipped with $s+1$ sections without common zeros $(g,\sigma_1,\ldots,\sigma_{s+1})$ the cobordism class of $(g,\sigma_1,\ldots,\sigma_{s+1},\tau_1,$ $\ldots,\tau_{r-s})$ for some sections $\tau_1,\ldots,\tau_{r-s}$. This, when considered as a map $\Prim^G_s(n,k)\to\Prim^G_r(n,k)$ is clearly the forgetful homomorphism $\varphi_s^r$ and so theorem \ref{thm:main} is proved.

\begin{rmk}
  \label{rmk:singpart}
  From the proof of lemma \ref{lemma:main} it is also apparent that the homomorphism $\Prim^G_r(n,k)\to\Imm^{(s+2)\zeta^G_{r-s-1}}(n-(s+1)(k+1),(s+2)(k+1))$ in theorem \ref{thm:main} assigns to the cobordism class of a prim $\Sigma^{1_r}$-map $f\colon M^n\to\R^{n+k}$ with immersion lift $g$, the cobordism class of $g|_{\ol\Sigma^{1_{s+1}}(f)}$ with its natural normal structure. Geometrically this means that the only obstruction to remove the singularities more complicated than $\Sigma^{1_s}$ of a prim $\Sigma^{1_r}$-map by a cobordism, is the cobordism class of its restriction to the closure of the $\Sigma^{1_{s+1}}$-locus.
\end{rmk}

\begin{rmk}
  The extreme cases $s=r-1$ and $r=\infty$ of theorem \ref{thm:main} are worthwhile to spell out; in fact these are the only cases we will need in the next section and they were (explicitly or implicitly) already obtained in \cite{ctrl}. Firstly, if $s=r-1$, we get the long exact sequence
  \begin{alignat*}2
    \ldots&\to\Prim^G_{r-1}(n,k)\to\Prim^G_r(n,k)\to\Imm^{(r+1)(\gamma^G_k\oplus\varepsilon^1)}(n-r(k+1),(r+1)(k+1))\to\\
          &\to\Prim^G_{r-1}(n-1,k)\to\ldots
  \end{alignat*}
  where $\varepsilon^1$ is the trivial line bundle and $\Imm^{(r+1)(\gamma^G_k\oplus\varepsilon^1)}(n-r(k+1),(r+1)(k+1))$ is naturally isomorphic to $\Imm^{(r+1)\gamma^G_k}(n-r(k+1),(r+1)k)$. This is the sequence coming from the analogue of the key fibration (see \cite{hosszu}) for prim cobordisms which was used in \cite{ctrl} to prove its main theorem. Secondly, if $r=\infty$, we get the long exact sequence
  \begin{alignat*}2
    \ldots&\to\Prim^G_s(n,k)\to\Imm^G(n,k+1)\to\Imm^{(s+2)\gamma^G_{k+1}}(n-(s+1)(k+1),(s+2)(k+1))\to\\
          &\to\Prim^G_s(n-1,k)\to\ldots
  \end{alignat*}
  which was shown to be rationally split in the case $G=\SO$ in \cite{ctrl} and was used to determine the ranks of the cobordism groups $\Prim^\SO_s(n,k)$.
\end{rmk}

\begin{rmk}
  For $s=0,r=\infty,k>0$ the exact sequence in theorem \ref{thm:main} is an exact sequence of Salomonsen \cite{salom}: the cobordism groups in this case are
  \begin{alignat*}2
    \Prim^G_0(n,k)&=\Imm^G(n,k),\\
    \Prim^G_\infty(n,k)&=\Imm^G(n,k+1),\\
    \Imm^{2\zeta^G_\infty}(n-k-1,2k+2)&=\Imm^{2\gamma^G_{k+1}}(n-k-1,2k+2)
  \end{alignat*}
  and the forgetful homomorphism $\varphi_0^\infty$ is the map induced by the embedding of a hyperplane $\R^{n+k}\subset\R^{n+k+1}$. Salomonsen's proof of his sequence is a special case of our proof here.
\end{rmk}

\section{Computations}
\label{sec:comp}

In the following we denote for any prime $p$ by $\CC_p$ the Serre class of Abelian groups of order a power of $p$ and for any number $m$ by $\CC_{\le m}$ the class of groups of order some combination of primes at most $m$. The symbol $\congc{}$ for a Serre class $\CC$ will denote isomorphism modulo $\CC$.

In \cite{nehez} the classifying spaces of cobordisms of codimension-$0$ oriented and unoriented and codimension-$1$ oriented prim maps were identified with familiar spaces, which yields the correspondences
\begin{alignat*}2
  &\Prim^\O_r(n,0)\cong\pi^s_{n+1}(\RP^{r+1})\congc2
  \begin{cases}
    0,&\text{if }r\text{ is even or }\infty,\\
    \pi^s(n-r),&\text{if }r\text{ is odd},
  \end{cases}\\
  &\Prim^\SO_r(n,0)\cong\pi^s(n)\oplus\pi^s(n-r),\\
  &\Prim^\SO_r(n,1)\cong\pi^s_{n+1}(\CP^{r+1})\congc{\le r+1}\displaystyle\bigoplus_{i=0}^r\pi^s(n-2i).
\end{alignat*}
Moreover, the groups $\Mor^G_r(n,1)$ for $G=\O,\SO$ were also computed modulo small torsion groups in \cite{szszt} and \cite{nulladik} respectively:
\begin{alignat*}2
  &\Mor^\O_r(n,1)\congc20,\\
  &\Mor^\SO_{2r+1}(n,1)\congc2\Mor^\SO_{2r}(n,1)\congc{\le2r+1}\displaystyle\bigoplus_{i=0}^r\pi^s(n-4i).
\end{alignat*}
We can now get an analogous description of codimension-$1$ unoriented prim cobordism groups:

\begin{prop}
  For all $n$ and $r$ we have
  $$\Prim^\O_{2r}(n,1)\congc2\Prim^\O_{2r-1}(n,1)\congc{\le2r}\displaystyle\bigoplus_{i=0}^{r-1}\pi^s(n-2-4i).$$
\end{prop}

\begin{prf}
  The exact sequence of theorem \ref{thm:main} in the case $G=\O,k=1,s=r-1$ is of the form
  $$\ldots\to\Prim^\O_{r-1}(n,1)\to\Prim^\O_r(n,1)\to\Imm^{(r+1)\gamma^\O_1}(n-2r,r+1)\to\Prim^\O_{r-1}(n-1,1)\to\ldots$$
  and since the Thom space $T(r+1)\gamma^\O_1$ is $\RP^\infty/\RP^r$, we have $\Imm^{(r+1)\gamma^\O_1}(n-2r,r+1)\cong\pi^s_{n-r+1}(\RP^\infty/\RP^r)$. Observe that if $r$ is even, then the homology of $\RP^\infty/\RP^r$ with $\Z_p$ coefficients vanishes for all odd primes $p$, hence by Serre's extension of the Whitehead theorem (see \cite{serre}), we have $\pi^s_{n-r+1}(\RP^\infty/\RP^r)\congc20$. This implies the isomorphism $\Prim^\O_r(n,1)\congc2$ $\Prim^\O_{r-1}(n,1)$ for all even $r$.

  The case of $r$ odd is strongly inspired by \cite{nulladik}. Note that if $r$ is odd, then the embedding $S^{r+1}\into\RP^\infty/\RP^r$ as a ``fibre'' of this Thom space is an isomorphism in homology with $\Z_p$ coefficients for all odd primes $p$, hence again by Serre's theorem, it is also an isomorphism in stable homotopy groups modulo $\CC_2$. Thus the long exact sequence above takes the form
  $$\ldots\to\Prim^\O_{r-1}(n,1)\to\Prim^\O_r(n,1)\to\pi^s(n-2(r+1))\to\Prim^\O_{r-1}(n-1,1)\to\ldots$$
  modulo $\CC_2$. We will prove that this sequence has an $m$-splitting $\pi^s(n-2(r+1))\to\Prim^\O_r(n,1)$, i.e. a homomorphism whose composition with the map $\Prim^\O_r(n,1)\to\pi^s(n-2(r+1))$ in the sequence is the multiplication by $m$ in $\pi^s(n-2(r+1))$, such that $m$ has no prime divisors larger than $r+1$. If such an $m$-splitting exists, then the long exact sequence splits modulo $\CC_{\le r+1}$, hence an induction on $r$ implies our statement.

  Now the following steps will conclude the proof:
  \begin{enumerate}[(i)]
    \item\label{11} For any $r$, if there is an immersion $M^{2r}\imto\R^{2r+2}$, where $M$ is oriented, such that the algebraic number of $(r+1)$-tuple points (i.e. the points counted with signs) is $m$, then there is an $m$-splitting $\sigma\colon\pi^s(n-2(r+1))\to\Prim^\O_r(n,1)$ of the above sequence.
    \item\label{12} For any $r$, the order $m(r)$ of the cokernel of the stable Hurewicz homomorphism
      $$h_{2r-2}\colon\pi^s_{2r-2}(\CP^\infty)\to H_{2r-2}(\CP^\infty)$$
      is such that there is an immersion $M^{2r}\imto\R^{2r+2}$ where $M$ is oriented and whose algebraic number of $(r+1)$-tuple points is $m(r)$.
    \item\label{13} For any $r$, the number $m(r)$ defined above has no prime divisors greater than $r+1$.
  \end{enumerate}

  All of these statements were essentially proved in \cite[parts 2 and 3]{nulladik}, so we only spell out what is different in our case. For the proof of (\ref{11}) the only thing missing in \cite{nulladik} is that the map constructed as a representative of the image $\sigma(x)$ for an element $x\in\pi^s(n-2(r+1))$ is actually a prim map; this is not spelled out there but it is apparent from its construction and remark \ref{rmk:prim}. The statement (\ref{12}) is only stated in \cite{nulladik} in the case where $r$ is even, however, its proof goes through for any $r$ by just changing the normal Pontryagin number $\la p_1(\nu_f)^{\frac r2},[M]\ra$ (for an immersion $f\colon M^{2r}\imto\R^{2r+2}$) to the normal Euler number $\la e(\nu_f)^r,[M]\ra$. Finally, the proof of (\ref{13}) in \cite{nulladik} goes through without change.
\end{prf}

\begin{rmk}
  The above proposition implies that the cobordism group $\Mor^\SO_r(n,1)$ is isomorphic to $\Prim^\O_{r+1}(n+2,1)$ modulo $\CC_{\le r+1}$ and the proof suggests an even stronger connection between them. It would be interesting to know if there is a geometric explanation for this.
\end{rmk}

This completes the computation of cobordism groups $\Prim^G_r(n,k)$ (modulo small torsion) when the codimension $k$ is at most $1$. Next we will consider the cases where $r=1$ and the codimension is large, more precisely when $n$ is either $2k+1$ or $2k+2$. Note that for $n\le2k$ we have $\Prim_1^G(n,k)\cong\Imm^G(n,k+1)$ since in this case a generic hyperplane projection of an immersion $M^n\imto\R^{n+k+1}$ only has $\Sigma^{1,0}$-type singularities and this also holds for cobordisms $W^{n+1}\imto\R^{n+k+1}\times[0,1]$ between them.

The cobordism groups $\Mor^G_1(n,k)$ were computed for $n=2k+1,2k+2$ in \cite{eszt} and \cite{2k+2} respectively (see also \cite[remark 10.11]{egzsor}) in terms of $\NN_n$ for $G=\O$ and $\Omega_n$ for $G=\SO$.
In our computation of their prim analogues we will use 
the group $\Imm^G(n,k+1)$ instead, which was computed modulo small torsion in \cite{immo} and \cite{immso} for $G=\O,\SO$ respectively:
\begin{alignat*}2
  &\Imm^\O(n,k+1)\congc2
    \begin{cases}
      \Omega_{n-k-1},&\text{if }k\text{ is odd},\\
      0,&\text{if }k\text{ is even},
    \end{cases}\\
  &\Imm^\SO(n,k+1)\congc{\le3}
    \begin{cases}
      P_{n,k+1}\oplus\Omega_{n-k-1},&\text{if }k\text{ is odd},\\
      P_{n,k+1},&\text{if }k\text{ is even},
    \end{cases}
\end{alignat*}
whenever $n<3k+3$, where $P_{n,k+1}$ is the subgroup of $\Omega_n$ consisting of cobordism classes $[M]$ whose each Pontryagin number $p_I[M]$ such that $p_I$ is divisible by a normal Pontryagin class $\ol p_i$ with $2i>k+1$, is zero.

\begin{lemma}
  \label{lemma:01}
  For all $k$ there is an exact sequence
  \begin{alignat*}2
    0&\to\Prim^G_1(2k+2,k)\to\Imm^G(2k+2,k+1)\to A\to\\
    &\to\Prim^G_1(2k+1,k)\to\Imm^G(2k+1,k+1)\to0.
  \end{alignat*}
  where $A\cong\Z$ if $G=\SO$ and $A\cong\Z_2$ if $G=\O$
  and the arrow $\Imm^G(2k+2,k+1)\to A$ is the homomorphism $\#\Sigma^{1_2}$ which assigns to the cobordism class of an immersion the algebraic number of $\Sigma^{1_2}$-points of its generic hyperplane projection (i.e. the number of points with their natural signs if $A\cong\Z$ and modulo $2$ if $A\cong\Z_2$).
\end{lemma}

\begin{prf}
  The exact sequence of theorem \ref{thm:main} with $r=\infty$ and $s=1$ contains the portion
  \begin{alignat*}2
    \ldots&\to\Imm^{3\gamma^G_{k+1}}(1,3k+3)\to\\
    \to\Prim^G_1(2k+2,k)\to\Imm^G(2k+2,k+1)&\to\Imm^{3\gamma^G_{k+1}}(0,3k+3)\to\\
    \to\Prim^G_1(2k+1,k)\to\Imm^G(2k+1,k+1)&\to0.
  \end{alignat*}
  Here we can use \cite[lemmas 1 and 2]{2k+2} to determine the cobordism group $\Imm^{3\gamma^G_{k+1}}(i,3k+3)\cong\pi^s_{3k+3+i}(T3\gamma^G_{k+1})$ for $i=0,1$; this yields that 
  for $i=0$ it is the group $A$ defined above and for $i=1$ it is trivial. The fact that $\Imm^G(2k+2,k+1)\to A$ is $\#\Sigma^{1_2}$, follows immediately from remark \ref{rmk:singpart}.
\end{prf}

\begin{prop}
  For all $k$ we have
  \begin{alignat*}2
    &\Prim^\O_1(2k+1,k)\cong\Imm^\O(2k+1,k+1),\\
    &\Prim^\O_1(2k+2,k)<\Imm^\O(2k+2,k+1)\text{ is of index }2.
  \end{alignat*}
\end{prop}

\begin{prf}  
  We use lemma \ref{lemma:01} with $G=\O$, by which it is sufficient to prove that the homomorphism $\#\Sigma^{1_2}\colon\Imm^\O(2k+2,k+1)\to\Z_2$ is non-trivial. In other words, we want to show that there is an immersion $f\colon M^{2k+2}\to\R^{3k+3}$ such that its hyperplane projection has an odd number of $\Sigma^{1_2}$-points. But by \cite{singprim} this number is (modulo $2$) the same as the number of triple points of $f$ and by \cite{3imm} for any $k$ there is an immersion $f$ as above with an odd number of triple points.
\end{prf}\medskip

In the following $\ell(n)$ will denote the largest number $l$ such that $3^l$ divides $\binom{3n}n$. For two groups $A,B$ we will use $A~?~B$ to denote an extension of $A$ by $B$, i.e. a group for which there is an exact sequence $0\to B\to A~?~B\to A\to0$.

\begin{prop}
  For all $k$ we have
  \begin{alignat*}2
    &\Prim^\SO_1(2k+1,k)\cong
      \begin{cases}
        \Imm^\SO(2k+1,k+1)~?~\Z_6,&\text{if }k=1,3,\\
        \Imm^\SO(2k+1,k+1)~?~\Z_{3^{\ell((k+1)/2)}},&\text{if }k>3\text{ is odd},\\
        \Imm^\SO(2k+1,k+1)~?~\Z,&\text{if }k\text{ is even},
      \end{cases}\\
    &\Prim^\SO_1(2k+2,k)\cong\Imm^\SO(2k+2,k+1).
  \end{alignat*}
\end{prop}

\begin{prf}
  We now use lemma \ref{lemma:01} with $G=\SO$; in this case we want to determine the homomorphism $\#\Sigma^{1_2}\colon\Imm^\SO(2k+2,k+1)\to\Z$.

  If $k$ is odd, then for any oriented immersion $f\colon M^{2k+2}\to\R^{3k+3}$ the triple points of $f$ (in the source $M$) have natural signs, and by \cite{singprim} we also have that their algebraic number is the same as $\#\Sigma^{1_2}[f]$. The minimal positive number of these triple points was determined by Eccles and Mitchell \cite{3immor} and Li \cite{3imm4}, using their notation it is $3\lambda_{\frac{k+1}2}$ where $\lambda_m$ is $2$ for $m=1,2$ and $3^{\ell(m)-1}$ for $m>2$ (here we need its triple since we are counting triple points in the source and not in the target). Thus the image of $\#\Sigma^{1_2}$ is $3\lambda_{\frac{k+1}2}\Z$ and our claim for $k$ odd follows.

  The case of $k$ even is simpler since again by \cite{singprim}, if $f\colon M^{2k+2}\to\R^{3k+3}$ is an oriented immersion where $k$ is even, then the oriented cobordism class of its $\Sigma^{1_2}$-points is of order $2$, i.e. it is trivial. Hence in this case the homomorphism $\#\Sigma^{1_2}$ is zero and the claim again follows.
\end{prf}

\end{document}